\documentclass{amsart}
\usepackage{amssymb,amsmath,amsthm}
\usepackage{mathrsfs}
\usepackage[all]{xy}
\begin{document}
\title{Talbot Workshop 2010 Talk 2: K-theory and Index Theory}
\author{Chris Kottke}
\email{ckottke@math.brown.edu}
\maketitle
\newtheorem*{prop}{Proposition}
\newtheorem*{thm}{Theorem}

\newcommand\pa{\partial}
\newcommand\id{\mathrm{Id}}
\newcommand\pt{\mathrm{pt}}
\newcommand\End{\mathrm{End}}
\newcommand\N{\mathbb{N}}
\newcommand\C{\mathbb{C}}
\newcommand\R{\mathbb{R}}
\newcommand\Z{\mathbb{Z}}

\newcommand\Sym{\mathrm{Sym}}

\newcommand\pair[1]{\left\langle #1\right\rangle}
\newcommand\set[1]{\left\{ #1 \right\}}
\newcommand\abs[1]{\left|#1\right|}
\newcommand\norm[1]{\left\|#1\right\|}
\newcommand\parens[1]{\left(#1\right)}
\newcommand\stackover[2]{\genfrac{}{}{-2pt}{2}{#1}{#2}}

\newcommand\wt[1]{\widetilde{{#1}}}

\newcommand\smallto{\rightarrow}
\renewcommand\to{\longrightarrow}
\renewcommand\mapsto{\longmapsto}

\newcommand\Diffeo{\mathrm{Diffeo}}
\newcommand\Ker{\mathrm{ker}}
\newcommand\Coker{\mathrm{coker}}
\newcommand\GL{\mathrm{GL}}
\newcommand\supp{\mathrm{supp}}
\newcommand\Hom{\mathrm{Hom}}
\newcommand\Iso{\mathrm{Iso}}
\newcommand\cV{\mathcal{V}}
\newcommand\cI{\mathcal{I}}
\newcommand\Id{\mathrm{Id}}
\newcommand\ind{\mathrm{ind}}
\newcommand\ch{\mathrm{ch}}
\newcommand\Td{\mathrm{Td}}
\newcommand\Thom{\mathfrak{T}}
\newcommand\cl{\mathrm{c}\ell}
\newcommand\Cl{\mathbb{C}\ell}
\newcommand\RCl{\mathrm{C}\ell}

\newcommand\Diff{\mathrm{Diff}}
\newcommand\sspan{\mathrm{span}}

\newcommand\smwedge{{\scriptstyle{\wedge}}}
\newcommand\swedge{{\textstyle{\bigwedge}}}
\newcommand\iprod{\lrcorner}
\newcommand\codim{\mathrm{codim}}
\newcommand\dimn{\mathrm{dim}}
\newcommand\mmod{\mathrm{\;mod\;}}
\newcommand\Fred{\mathrm{Fred}}
\newcommand\Spin{\mathrm{Spin}}
\newcommand\SO{\mathrm{SO}}
\newcommand\rmO{\mathrm{O}}
\newcommand\bbS{\mathbb{S}}
\newcommand\dirD{{/}\!\!\!\!D}
\newcommand\bbE{\mathbb{E}}
\newcommand\cM{\mathcal{M}}
\newcommand\hcM{\widehat{\mathcal{M}}}

From the point of view of an analyst, one of the most delightful things about complex K-theory is that it has a nice realization by analytical objects, namely
(pseudo)differential operators and their Fredholm indices.  This connection allows quite a bit of interesting information to flow both ways: from analysis to topology
and vice versa.  

This talk will try and give a sketch of this picture, and consists of three parts or themes.  The first is ``the Gysin map as the index,'' describing the 
families index theorem of Atiyah and Singer, and how the pushforward along a fibration in K-theory can be realized as the index of a family of operators.  The 
second is ``spin$^c$ as an orientation,'' in which I discuss Clifford algebras, spin and spin$^c$ structures, Dirac operators and the analytic realization
of the Thom isomorphism for complex K-theory.  Finally (I did not have time to get to this part in the actual Talbot talk) I will discuss the constructions
leading to higher index maps (i.e. $K^1$ instead of $K^0$; of course this is more interesting for real K-theory than for complex K-theory), namely Clifford-linear
differential and Fredholm operators.  The best reference for almost everything in this talk is the wonderful book \cite{lawson1989spin} by Lawson and Michelson.  
I cannot recommend this book highly enough.  I will also try and give references to original sources (essentially all of which involve Atiyah as an author).

\section{Some notation and facts} \label{S:some}
First let us get down some notation and facts about K-theory that will be of use in the following.  Let $V \to X$ be a (not necessarily complex) vector bundle.  
There are many ways of constructing the {\bf Thom space} of $V$, which will be denoted $X^V$:
$$
	X^V := DV/SV = \overline{V}/\pa V = \mathbb{P}(V\oplus 1)/\infty.
$$
The first space denotes the unit disk bundle of $V$ (with respect to some choice of metric) quotiented out by the unit sphere bundle; the second denotes the
radial compactification of $V$ quotiented out by its boundary; the third denotes the projective bundle\footnote{Note that this constructs the fiberwise one point 
compactification of $V$.} associated to $V\oplus 1$ (where $1$ is the trivial complex line bundle), modulo the section at infinity.  Indeed we could take 
any compactification of $V$ modulo the points added at infinity.

In any cohomology theory, the reduced cohomology of $X^V$ is of interest.  In index theory, however, we prefer to think of the K-theory of $X^V$ as 
{\bf compactly supported K-theory}\footnote{Though we shall only need it for vector bundles, by compactly supported cohomology of any space $M$ 
can be thought of as the relative cohomology of $(M^+, \infty)$ where $M^+$ denotes the (one-point or otherwise) compactification of $M$.  This is consistent
if we agree to take $M^+ = M\sqcup \set{pt}$ for compact $M$.} of the space $V$:
$$
	K^\ast_c(V) := \wt K^\ast(X^V) = K^\ast(\overline V, \pa \overline V).
$$

There is a convenient representation of relative even K-theory of a pair $(X, A)$, where $A \subset X$ is a nice enough subset, known as the {\bf difference bundle
construction}:
$$
	K^0(X, A) = \set{E,F,\sigma \;;\; \sigma_{|A} : E \stackrel{\cong}{\to} F}/\sim,
$$
where $E,F \to X$ are vector bundles and $\sigma : E \to F$ is a bundle map covering the identity on $X$, which restricts to an isomorphism over $A$.  The equivalence
relations amount to stabilization and homotopy.  Intuitively this should be clear; if $E$ and $F$ are isomorphic over $A$, then, $[E] - [F]$ should be trivial
in K-theory when restricted to $A$.

Unpacking this in the case of compactly supported K-theory for $V$, we conclude that we can represent $K^0_c(V)$ by
$$
	K^0_c(V) = \set{[\pi^\ast E, \pi^\ast F, \sigma ]\;;\; \sigma_{V\setminus 0} : \pi^\ast E \stackrel{\cong}{\to} \pi^\ast F},
$$
where $\pi : V\to X$ is the projection and $\mathbf{0}$ denotes the zero section.  Indeed, the pair $(\overline V, \pa \overline{V})$ is homotopy equivalent to 
$(V, V \setminus \mathbf{0})$, and by contractibility of the fibers of $\overline V$, any vector bundles over $V$ are homotopic to ones pulled up from the base, 
i.e.\ of the form $\pi^\ast E$.

Let $V \to X$ now be a complex vector bundle.  The {\bf Thom isomorphism} in K-theory states that $V$ has a K-theory orientation, so that 
$\wt K^\ast (X^V) = K^\ast (X)$.  Specifically, $\wt K^\ast(X^V) = K^\ast_c(V)$ is a freely generated, rank one module over $K^\ast(X)$, and in the representation 
of compactly supported K-theory discussed above, the generator, or {\bf Thom class}\footnote{We'll discuss orientation classes more generally in section \ref{S:gysin}, 
and we'll interpret the Thom class in terms of spin$^c$ structures in section \ref{S:structures}.}
$\mu \in K^0_c(V)$ has the following nice description.
\begin{prop}
The Thom class $\mu \in K^0_c(V)$ for a complex vector bundle $V \to X$ can be represented as the element
$$
	\mu = [\pi^\ast \swedge^\mathrm{even}V, \pi^\ast \swedge^\mathrm{odd} V, \cl] \in K^0_c(V), \qquad \cl(v)\cdot = v\smwedge \cdot - v^\ast \iprod \cdot,
$$
where the isomorphism off $0$ is given by $\sigma(v) = \cl(v) = v \smwedge \cdot - v^\ast \iprod \cdot$, the first term denoting exterior product with $v$ and
the second denoting the contraction with $v^\ast$ (equivalently, the inner product with $v$), with respect to any choice of metric.  
\end{prop}
The isomorphism $\cl(v)$ is an example of Clifford multiplication, about which we will have much more to say in section \ref{S:clifford}.

Finally, a bit about Fredholm operators.  Let $H$ be a separable, infinite dimensional Hilbert space, and recall that a bounded linear operator $P$ is
{\bf Fredholm} if it is invertible modulo {\bf compact operators}, which in turn are those operators in the norm closure of the finite rank operators.
Thus $P$ is Fredholm iff there exists an operator $Q$ such that $PQ - \id$ and $QP - \id$ are compact.  One upshot of this is that
$$
	\Ker(P) \text{ and } \Coker(P) = \Ker(P^\ast) \quad \text{are finite dimensional.}
$$
The relationship between Fredholm operators and K-theory starts with the observation of Atiyah \cite{atiyah1967ktheory} that the space of Fredholm operators 
on $H$ classifies $K^0$:
\begin{prop} [Atiyah]
$$
	[X, \Fred(H)] = K^0(X),
$$
where the left hand side denotes homotopy classes of maps $X \to \Fred(H)$, the latter given the operator topology.\footnote{The precise topology one should take on 
$\Fred(H)$ becomes a little difficult in twisted K-theory, but (I guess!) not here.}  
\end{prop}
Morally, the idea is to take a map $P : X \to \Fred(H)$, and examine the
vector bundles $\Ker(P)$ and $\Coker(P)$, whose fibers at a point $x \in X$ are the finite dimensional vector spaces $\Ker(P(x))$ and $\Coker(P(x))$, respectively.
Of course this is a bit of a lie, since the ranks of these bundles will generally jump around as $x$ varies; nevertheless, it is possible to stabilize the
situation and see that the class
$$
	[\Ker(P)] - [\Coker(P)] \in K^0(X)
$$
is well-defined.

\section{Differential operators and families} \label{S:differential}
One of the most important sources of such maps $X \to \Fred(H)$ are families of differential operators on $X$.  Let's start with differential operators
themselves.  A working definition of the {\bf differential operators of order $k$}, $\Diff^k(X; E,F) : C^\infty(X; E)\to C^\infty(X; F)$, where $E$ and $F$ are
vector bundles over $X$ is the following local definition.
$$
	\Diff^k(X; E,F) \ni P \stackrel{\text{locally}}{=} \sum_{\abs{\alpha} \leq k} a_\alpha(x) \pa_x^\alpha, \quad a_\alpha(x) \in \Hom(E_x,F_x),
$$
where we're employing multi-index notation: $\alpha = (\alpha_1, \ldots,\alpha_n) \in \N^n$, $\abs{\alpha} = \sum_i \alpha_i$, 
$\pa_x^\alpha = \pa_{x_1}^{\alpha_1}\cdots \pa_{x_n}^{\alpha_n}$, where $x = (x_1,\ldots,x_n)$ are local coordinates on $X$.  

This local expression for $P$ does not transform well under changes of coordinates; however, the highest order terms (those with $\abs{\alpha} = k$) {\em do} 
behave well.  If we consider $\pa_x^\alpha \in \Sym^{\abs{\alpha}} T_x X$, we can view it as a monomial map $T_x^\ast X \to \R$ of order $\abs{\alpha}$.  
If $x = (x_1,\ldots,x_n)$ are coordinates on $X$ inducing coordinates $(x,\xi) = (x_1,\ldots,x_n,\xi_1,\ldots,\xi_n)$ on $T^\ast X$, the monomial obtained is just
$$
	\pa_x^\alpha = \xi^\alpha = \xi_1^{\alpha_1}\cdots \xi_n^{\alpha_n}.
$$
Summing up all the terms 
of order $k$ gives us a homogeneous polynomial of order $k$, which because of the $a_\alpha$ term is a homogeneous polynomial on $T_x^\ast X$ valued in 
$\Hom(E_x,F_x)$.  The claim is that this {\bf principal symbol} 
$$
	\sigma(P)(x,\xi) = \sum_{\abs{\alpha} = k} a_\alpha(x) \xi^\alpha \in C^\infty(T^\ast X; \Hom(\pi^\ast E,\pi^\ast F))
$$
is well-defined.  

An operator is {\bf elliptic} if its principal symbol is invertible away from the zero section $\mathbf{0} \in T^\ast X$.  The canonical example of an elliptic operator
is $\Delta$, the Laplacian (on functions, say), a second order operator whose principal symbol is $\sigma(\Delta)(\xi) = \abs{\xi}^2$, where $\xi \in T^\ast X$ and the 
norm comes from a Riemannian metric.  The canonical non-example on $\R\times X$ is $\Box = \pa_t^2 - \Delta$, the D'Alembertian or wave-operator, whose 
principal symbol is $\sigma(\Box) = \tau^2 - \abs{\sigma}^2$, where $(\tau,\sigma) \in T^\ast \R \times T^\ast X$, which vanishes on the (light) cones 
$\set{\tau = \pm\abs{\xi}}$.

The reader whose was paying particularly close attention earlier will note that the symbol of an elliptic differential operator is just the right kind of object to 
represent an element in the compactly supported K-theory\footnote{In fact, any element of $K_c^0(T^\ast X)$ can be represented as the symbol of an elliptic 
{\em pseudodifferential} operator, though we shall not discuss these here.} of $T^\ast X$ since it is 
invertible away from $\mathbf{0} \subset T^\ast X$:
$$
P \text{ elliptic} \implies [\pi^\ast E, \pi^\ast F, \sigma(P)] \in K_c^0(T^\ast X).
$$

For our purposes, the other important feature of an elliptic operator $P$ on a compact manifold is that it extends to a Fredholm 
operator $P: L^2(X; E) \to L^2(X; F)$.  Actually, this is a bit of a lie, since if $k > 0$, $P \in \Diff^k(X; E, F)$ is unbounded on $L^2(X; E)$, and we should 
really consider it acting on its maximal domain in $L^2(X; E)$, which is the Sobolev space $H^k(X; E)$ which itself has a natural Hilbert space structure.  However, 
since the order of operators is immaterial 
as far as index theory is concerned, we will completely ignore this issue for the rest of
this note, pretending all operators in sight are of order zero\footnote{In fact it is always possible to compose $P$ with an invertible {\em pseudodifferential} 
operator (of order $-k$) so that the composite has order zero, without altering the index of $P$.  One such choice is $(1+\Delta)^{-k/2}$, another
is $(1 + P^\ast P)^{-1/2}$.}, which act boundedly 
on $L^2$.  

Now let $X \to Z$ be a fibration of compact manifolds with fibers $X_z \cong Y$.  A {\bf family of differential operators} with respect to $X \to Z$, is just 
a set of differential
operators on (vector bundles over) the fibers $X_z$, parametrized smoothly by the base $Z$.  For a formal definition, take the principal $\Diffeo(Y)$ bundle 
$\mathcal{P} \to Z$ such that $X = \mathcal{P} \times_{\Diffeo(Y)} Y$; then the differential operator families of order $k$ are obtained as\footnote{I will be sloppy 
about distinguishing between families vector bundles on the fibers and vector bundles on the total space $X$.  In fact they are the same.}
$$
	\Diff^k(X/Z; E_1,E_2) = \mathcal{P} \times_{\Diffeo(Y)} \Diff^k(Y; E_1, E_2).
$$

As before, there is a principal symbol map
$$
	\Diff^k(X/Z; E_1,E_2) \ni P \mapsto \sigma(P) \in C^\infty(T^\ast(X/Z); \Hom(\pi^\ast E_1, \pi^\ast E_2)),
$$
where $T^\ast(X/Z)$ denotes the vertical (a.k.a.\ fiber) cotangent bundle.  Once again $P$ is {\bf elliptic} if $\sigma(P)$ is invertible away from the zero section; if 
this is the case, $P$ extends to a family of Fredholm operators on the Hilbert space bundles
$$
	\mathcal{H}_i \to Z \quad i = 1,2 \quad \text{with fiber $(\mathcal{H}_i)_z = L^2(X_z; E_i)$.}
$$
By Kuiper's theorem that the unitary group of an infinite dimensional Hilbert space is contractible, the bundles $\mathcal{H}_i$ are trivializable, so trivializing and 
identifying $\mathcal{H}_1$ and $\mathcal{H}_2$ (all separable, infinite dimensional Hilbert spaces are isomorphic), we obtain a map
$$
	P : Z \to \Fred(H), \quad H \cong L^2(Y)
$$
which must therefore have an index in the even K-theory of $Z$:
$$
	\Diff^\ast(X/Z; E, F) \ni P \text{ elliptic} \implies \ind(P) = [\Ker(P)] - [\Coker(P)] \in K^0(Z).
$$

Just as in the case of the single operator, the principal symbol of the family $P$ represents a class in compactly supported K-theory of $T^\ast (X/Z)$:
$$
	[\pi^\ast E, \pi^\ast F, \sigma(P)] \in K^0_c(T^\ast(X/Z)).
$$
We will come back to the relationship between these two objects in a moment; for now you should think of the index as an assignment which maps
$[\pi^\ast E, \pi^\ast F, \sigma(P)] \in K^0_c(T^\ast(X/Z))$ to $\ind(P) = [\Ker(P)] - [\Coker(P)] \in K^0(Z)$.  This is well-defined since any two elliptic operators
with the same principal symbol are homotopic through elliptic (hence Fredholm) operators, and since the index is homotopy invariant, any two choices of operators
$P,P'$ with the same symbol $\sigma(P) = \sigma(P')$ will have the same index in $K^0(Z)$.

\section{Gysin maps}\label{S:gysin}
In the first Talbot talk, Jesse Wolfson discussed the Gysin map in K-theory associated to an embedding.  We will need a similar kind of Gysin map associated to
fibrations.  Let $p : X \to Z$ be a smooth fibration with $Z$ compact but not necessarily having compact fibers.  By the theorem of Whitney, we can embed 
any manifold into $\R^N$ for sufficiently large $N$, and since $Z$ is compact, this can be done fiberwise to obtain an embedding of fibrations from $X$
into a trivial fibration:
$$
\xymatrix{
	X \ar@{^{(}->}[rr] \ar[dr]^{p} && Z \times \R^N\ar[dl]_{\mathrm{pr}_1} \\ & Z &
}
$$
Let $\nu \to X$ denote the normal bundle to $X$ with respect to this embedding; by the collar neighborhood theorem it is isomorphic to an open neighborhood of 
$X$ in $Z\times \R^N$.  The open embedding $i : \nu \hookrightarrow Z\times \R^N$ induces a wrong way map
$$
	 \tilde i : \Sigma^N Z \to X^\nu
$$
by adding points at infinity and considering the quotient map $Z\times \R^N / \infty \to \nu / \infty$.  Thus we obtain a {\bf Gysin (a.k.a.\ ``wrong way,'' 
``umkehr,'' ``pushforward,'' ``shriek'') map}
$$
	\tilde i^\ast : \wt h^\ast(X^\nu) \to \wt h^\ast(\Sigma^N Z) = h^{\ast - N}(Z)
$$
in any generalized cohomology theory $h^\ast(\cdot)$.  If the fibration is oriented (in a sense defined below), we will actually obtain a map 
from the cohomology of $X$ to that of $Z$.

Let $V \to X$ be a vector bundle.  We say that $V$ has an {\bf orientation for the cohomology theory $h^\ast$} (which we are assuming is multiplicative) if there 
is a global (Thom) class $\mu \in \wt h^n\parens{X^V}$ which restricts to the multiplicative unit $\mu_x \in \wt h^n(X_x^V) = \wt h^n(S^n) = \wt h^0(\pt)$ 
(here $n$ is the rank of the vector bundle); if such a class exists, we have a Thom isomorphism 
$$
	h_c^\ast(X) \stackrel{\cong}{\to} \wt h^{\ast + n}\parens{X^V}.
$$
Specifically, $\wt h^\ast(X^V)$ is a freely generated module over $h_c^\ast(X)$ with generator $\mu$.

We say a fibration $X \to Z$ is {\bf oriented} with respect to the cohomology theory if $T(X/Z) \to X$ has an $h^\ast$ orientation.  Indeed, if this is the case,
the orientation on $T(X/Z) \to X$ induces\footnote{It is a theorem that if $\alpha$ and $\beta$ are vector bundles over $X$, than orientability of any
two of $\alpha, \beta, \alpha\oplus\beta$ implies orientability of the third.} one on $\nu \to X$, and we obtain the {\bf Gysin map associated to an oriented 
fibration}
$$
	p_! : h_c^\ast(X) \to h^{\ast - n}(Z)
$$
via the composition $h_c^\ast(X) \stackrel{\cong}{\to} \wt h^{\ast + (N-n)}(X^\nu) \stackrel{\tilde i^\ast}{\to} h^{\ast - n}(Z)$.  Note in particular that the 
degree shifts by $N$ cancel; indeed the Gysin map is completely independent of the choice of embedding $X \to Z\times \R^N$.

\section{The Gysin map as the index} \label{S:the}
Now let us return to families of elliptic differential operators.  Given $P \in \Diff^k(X/Z; E,F)$, we have the element
$[\pi^\ast E, \pi^\ast F, \sigma(P)] \in K^0_c(T^\ast(X/Z))$, which maps to the index $\ind(P) = [\Ker(P)] - [\Coker(P)] \in K^0(Z)$.  The famous index theorem
of Atiyah and Singer \cite{atiyah1968index_I_III} \cite{atiyah1971index_IV} can now be stated quite simply.
\begin{thm} [Atiyah-Singer]
The index map 
$$
	\ind : [\pi^\ast E,\pi^\ast F, \sigma(P)] \in K^0_c(T^\ast(X/Z)) \to [\Ker(P)] - [\Coker(P)] \in K^0(Z)
$$
coincides with the Gysin map $K_c^0(T^\ast(X/Z)) \to K^0(Z)$ associated to the oriented fibration\footnote{We'll see below why this fibration is canonically oriented.}
$$
	p : T^\ast (X/Z) \to Z.
$$
In short,
$$
	\ind = p_! : K^0_c(T^\ast(X/Z)) \to K^0(Z).
$$
\end{thm}
In particular, we can recover the integer index $\ind(P) = \dimn\; \Ker(P) - \dimn\;\Coker(P) \in \Z$ of a single operator on a compact manifold $X$ from the 
case $Z = \pt$; from the unique map $X \to \pt$, we get an oriented fibration $T^\ast X \to \pt$ and a Gysin map $p_! : K^0_c(T^\ast X) \to K^0(\pt) = \Z$.  

Let us unpack this a bit.  We have the fibration $T^\ast (X/Z) \to Z$ which factors as $T^\ast(X/Z) \to X\to Z$.  As noted above there is always an embedding
of $X$ into a trivial Euclidean bundle $Z\times \R^N \to Z$.  This induces an embedding $T^\ast (X/Z) \hookrightarrow Z\times T^\ast \R^N = Z\times \R^{2N}$, 
so we have the following situation
$$
\xymatrix{
	T^\ast(X/Z) \ar[d] \ar@{^{(}->}[r] & Z\times \R^{2N} \ar[d] \\ X \ar[d] \ar@{^{(}->}[r] & Z\times \R^N \ar[dl] \\ Z
}
$$

The point is that, because the embedding $T^\ast(X/Z) \hookrightarrow Z\times\R^{2N}$ comes from an embedding of $X$, the normal bundle $\nu \to T^\ast(X/Z)$
carries a canonical complex structure.  Indeed, if we denote by $NX \to X$ the normal bundle of $X$ with respect to $X \hookrightarrow Z\times \R^N$, then
$\nu$ is isomorphic to two copies of $NX$:
$$
	\nu \cong NX\oplus NX \cong NX \otimes \C,
$$
one copy representing the normal to the base $X$, and the other copy representing the normal to the fiber $T^\ast(X/Z)_x, x \in X$.  Thus we conclude that, whether
or not the fibration $X \to Z$ has a K-theory orientation, $T^\ast(X/Z) \to Z$ {\em always} has a K-theory orientation, so we have a Gysin map
$$
	p_! : K^0_c(T^\ast(X/Z)) \to K^0(Z),
$$
which coincides with the index map on elements in $K^0_c(T^\ast(X/Z))$ which represent symbols of elliptic differential operators\footnote{As remarked in 
previous footnotes, it is desirable to broaden one's focus to include pseudodifferential operators, for then every element of $K^0_c(T^\ast(X/Z))$ can
be represented as the symbol of an elliptic pseudodifferential operator which extends to a Fredholm operator, so the index map extends to all of 
$K^0_c(T^\ast(X/Z))$ and is equal to the Gysin map.}.  Note that there is no degree shift since $T^\ast(X/Z)$ has even dimensional fibers over $Z$ and
K-theory is 2-periodic.

In the next sections we shall discuss the conditions necessary for $X \to Z$ to have a K-theory orientation, and how to realize the Gysin map 
$K^0(X) \to K^0(Z)$ in terms of elliptic differential operators and their indices; this will involve a digression through Clifford algebras, spin groups,
spin structures and Dirac operators.  Later we will see how to deal with objects in odd K-theory.

\section{Clifford algebras} \label{S:clifford}
Let $(V,q)$ be a finite dimensional vector space over $\R$ with $q$ a non-degenerate quadratic form.  The (real) {\bf Clifford algebra} $\RCl(V,q)$ is 
the universal object with respect to maps $f : V \to A$, where $A$ is an associative algebra with unit, satisfying $f(v)\cdot f(v) = - q(v)1$.  
It can be constructed as a quotient of the tensor algebra:
$$
	\RCl(V,q) = \bigoplus_{n=0}^\infty V^{\otimes n} / \mathcal I, \quad \mathcal I = \pair{v\otimes v + q(v)1}
$$
where $\mathcal I$ is the ideal generated by all elements of the form $v\otimes v + q(v)1$.

As a vector space (but not as an algebra unless $q \equiv 0$!) $\RCl(V,q)$ is isomorphic to the exterior algebra $\bigoplus_{n=0}^\dimn(V) \bigwedge^n V$; 
in particular if $\set{e_i}$ is a basis of $V$, then $\set{e_{i_1}\cdots e_{i_k} \;;\; i_1 < \cdots < i_k}$ form a basis for $\RCl(V,q)$, which under multiplication 
are subject to the relation
$$
	e_i e_j = - e_j e_i - 2 q(e_i,e_j)
$$
where we denote also by $q$ the bilinear form associated to $q$.  Computations are easiest when $\set{e_i}$ is an orthonormal basis, whence the multiplication
simplifies to the rules $e_i e_j = -e_j e_i, i \neq j$ and $e_i^2 = -1$.

In fact $\RCl(V,q)$ is a $\Z_2$-graded algebra.  Indeed, if we let $\RCl^0(V,q)$ and $\RCl^1(V,q)$ be the images of $\bigwedge^\mathrm{even} V$ and 
$\bigwedge^\mathrm{odd} V$, respectively, under the vector space isomorphism with the exterior algebra, it is easy to check that 
$$
	\RCl(V,q) = \RCl^0(V,q)\oplus\RCl^1(V,q) \quad \text{and} \quad \RCl^i(V,q)\cdot \RCl^j(V,q) \subset \RCl^{(i+j)\mmod 2}(V,q).
$$
Alternatively, the involution $\alpha : V \to V : v \mapsto -v$ extends multiplicatively to an involution on all of $\RCl(V,q)$:
$$
	\alpha : \RCl(V,q) \to \RCl(V,q), \quad \alpha^2 = \id, \quad \alpha : V \ni v \mapsto -v \in V
$$
and we can define $\RCl^0(V,q)$ and $\RCl^1(V,q)$ as the positive and negative eigenspaces of $\alpha$, respectively.  Note in particular that 
$V \subset \RCl^1(V,q)$ as a vector subspace.  

We can also form the {\bf complex Clifford algebra} $\Cl(V,q)$ by tensoring up with $\C$:
$$
	\Cl(V,q) := \RCl(V\otimes \C,q_{\C}) \cong \RCl(V,q)\otimes \C.
$$
Complex Clifford algebras are those of primary importance for this talk, since it concerns (mostly) complex K-theory.  There is a parallel relationship
between real Clifford algebras and real K-theory.

We will denote the Clifford algebra of Euclidean $n$-space by
$$
	\RCl_n := \RCl(\R^n, \pair{\cdot,\cdot})
$$
and call it the {\bf (real) Clifford algebra of dimension $n$}.  Similarly, we will denote the {\bf complex Clifford algebra of dimension $n$} by
$$
	\Cl_n := \Cl(\R^n,\pair{\cdot,\cdot}) = \RCl(\C^n, \pair{\cdot,\cdot}).
$$

It is a rather nice fact that complex Clifford algebras are isomorphic to matrix algebras\footnote{In fact real Clifford algebras are also isomorphic to 
matrix algebras over $\R$, $\C$ or $\mathbb{H}$, or a direct sum of two such, with an 8-periodic pattern related to Bott periodicity in the real setting.}.
\begin{prop}
$$
	\Cl_{2n} \cong M(2^n, \C) \quad \text{and} \quad \Cl_{2n+1} \cong M(2^n,\C)\oplus M(2^n,\C)
$$
where $M(k, \C)$ denotes the algebra of $k\times k$ complex matrices.  
\end{prop}
As a consequence of this, the representations of $\Cl_n$ are easy to classify:
$\Cl_{2n}$ has a unique irreducible representation of (complex) dimension $2^n$ given by the obvious action of $M(2^n,\C)$ on $\C^{2^n}$; and
$\Cl_{2n+1}$ has two distinct irreps of dimension $2^n$ corresponding to action of one or the other of the factors of $M(2^n,\C)$.

The last tidbit we shall need is the algebra isomorphism 
$$
	\Cl_{n-1} \cong \Cl^0_n 
$$ 
(note that $\Cl^0_n$ is a subalgebra of $\Cl_n$).  This is obtained by considering the map $f : \R^{n-1} \to \Cl^0_n$ given on basis vectors by
$$
	\R^{n-1} \ni e_i \mapsto f(e_i) = e_i\,e_{n} \in \Cl^0_n.
$$
This satisfies $f(v)\cdot f(v) = v\,e_n\,v\,e_n = - q(v)1$ and thus generates an algebra map $\Cl_{n-1} \to \Cl^0_n$ by the universal property, which is
easily seen to be bijective.

This isomorphism leads to an equivalence between graded $\Cl_n$ modules and ungraded $\Cl_{n-1}$ modules.  In the one direction, if $M = M^0\oplus M^1$ is a graded
module over $\Cl_n$, then $M^0$ and $M^1$ are (possibly inequivalent!) modules over $\Cl_{n-1} \cong \Cl^0_n$.  In the other direction, 
given a $\Cl_{n-1}$-module $M$, the we can form $M\otimes_{\Cl^0_n} \Cl_n$, which is a graded module over $\Cl_n$.  

Putting this fact together with the classification of irreps above, we see that, for even Clifford algebras, there is a unique irreducible $\Cl_{2n}$ module 
$$
	\bbS_{2n} = \bbS^+_{2n}\oplus\bbS^-_{2n}
$$
which splits as the two inequivalent irreps of $\Cl_{2n}^0 \cong \Cl_{2n - 1}$.  Hence it has two {\bf inequivalent} gradings, either 
$\bbS^0\oplus\bbS^1 = \bbS^+\oplus\bbS^-$ or $\bbS^0\oplus\bbS^1 = \bbS^-\oplus\bbS^+$.  On the other hand, for odd Clifford algebras, there is a unique 
{\em graded} $\Cl_{2n + 1}$ module
$$
	\bbS_{2n + 1} = \bbS^+_{2n + 1} \oplus \bbS^-_{2n + 1}
$$
since both $\bbS^+_{2n + 1}$ and $\bbS^-_{2n + 1}$ must be equivalent to the unique irrep of $\Cl^0_{2n + 1} \cong \Cl_{2n}$.

\section{Spin and Spin$^c$ groups} \label{S:spin}
Given $(V,q)$, the group $\Spin(V,q)$ is the universal cover of the special orthogonal group $\SO(V,q)$.  We can find it inside the Clifford algebra of $V$ as follows.
Let $\RCl(V,q)^\times$ be the group of units inside $\RCl(V,q)$.  This group acts on $\RCl(V,q)$ by a twisted conjugation:
$$
	\RCl(V,q)^\times \times \RCl(V,q) \ni (x,v) \mapsto x\,v\,\alpha(x)^{-1}
$$
where $\alpha_{|\RCl^i(V,q)} = (-1)^i\id$ is the involution from earlier.  The {\bf Clifford group} $\Gamma \subset \RCl(V,q)^\times$ is the subgroup which fixes
the subspace $V \subset \RCl(V,q)$; it also preserves the quadratic form $q$ and hence maps to the orthogonal group $\rmO(V,q)$ with kernel $\R^\times$:
$$
	1 \to \R^\times \to \Gamma \to \rmO(V,q) \to 1.
$$
Up to a scalar factor, there is a natural choice of multiplicative norm $\abs{\cdot} : \Gamma \to \R^\times$, and the {\bf Spin group of $(V,q)$}, $\Spin(V,q)$, 
is defined to be the 
subgroup of norm 1 elements covering\footnote{The subgroup of norm 1 elements covering $\rmO(V,q)$ is called $\mathrm{Pin}(V,q)$, a joke which is apparently due to Serre.} 
the special orthogonal group $\SO(V,q)$:
$$
	\Spin(V,q) := \set{u \in \Gamma \;;\; \abs{u} = 1, u \text{ maps to } \SO(V,q)} \subset \Gamma.
$$
Alternatively, it can be defined as the subgroup of $\RCl(V,q)^\times$ generated by finite products of the form $v_1\cdots v_{2n}$, $v_i \in V$, $q(v_i) = 1$ with an 
even number of factors.  We have the exact sequence
$$
	1 \to \set{\pm 1} \to \Spin(V,q) \to \SO(V,q) \to 1
$$
and $\Spin(V,q)$ is compact if $q$ has positive signature.  It also lies in the 0-graded component of $\RCl(V,q)$:
$$
	\Spin(V,q) \subset \RCl^0(V,q)
$$
The spin group of $(\R^n,\pair{\cdot,\cdot})$ will simply be called {\bf the spin group of dimension $n$}, and denoted
$$
	\Spin_n := \Spin(\R^n, \pair{\cdot,\cdot}).
$$

From now on, we focus on the even dimensional case.  Since $\Spin_{2n} \subset \RCl_{2n}$, we have a complex representation coming from the 
irreducible $\Cl_{2n} = \RCl_{2n}\otimes \C$ module $\bbS_{2n}$.  In fact, since $\Spin_{2n} \subset \RCl^0_{2n}$, this splits as two inequivalent, irreducible
{\bf half spin representations}
$$
	\bbS_{2n} = \bbS^+_{2n} \oplus \bbS^-_{2n}.
$$
We denote these fundamental representations by
$$
	\rho_{1/2}^\pm : \Spin_{2n} \to \GL(\bbS^\pm_{2n}).
$$
The representation
$$
	\rho := \rho^+_{1/2}\oplus \rho^-_{1/2} : \Spin_{2n} \to \GL(\bbS_{2n}) = \GL(\bbS^+_{2n}\oplus\bbS^-_{2n})
$$
is called the {\bf fundamental spin representation} and vectors in $\bbS_{2n}$ are called {\bf spinors}.  We'll see the importance of spinors in defining K-theory
orientation classes in section \ref{S:structures}.

As we are primarily interested in complex representations, there is another group inside $\Cl_n$ to discuss.  {\bf The  Spin$^c$ group associated to $(V,q)$}
is the quotient
$$
	\Spin^c(V,q) = \Spin(V,q) \times_{\Z_2} \mathrm{U}_1
$$
where the $\Z_2$ is generated by the element $(-1,-1) \in \Spin(V,q)\times \mathrm{U}_1$.  We have the exact sequence
$$
	1 \to \set{\pm 1} \to \Spin^c(V,q) \to \SO(V,q)\times \mathrm{U}_1 \to 1
$$
As with the spin group, $\Spin^c(V,q)$ sits inside $\Cl(V,q)$,
$$
	\Spin^c(V,q) \subset \Cl^0(V,q) \subset \Cl(V,q) = \RCl(V,q) \otimes \C
$$
and, for the canonical spin$^c$ groups of even dimension, $\Spin^c_{2n} := \Spin^c(\R^{2n}, \pair{\cdot,\cdot})$, we have a {\bf fundamental spin$^c$ representation}
$$
	\rho := \rho^+_{1/2}\oplus \rho^-_{1/2} : \Spin^c_{2n} \to \GL(\bbS_{2n}) = \GL(\bbS^+_{2n}\oplus\bbS^-_{2n})
$$
on the {\bf spinors} $\bbS_{2n}$, where again, $\bbS_{2n}$ is the unique irreducible $\Cl_{2n}$ module.

\section{Spin$^{(c)}$ structures} \label{S:structures}
Let us now transfer this to a manifold setting.  We'll define spin and spin$^c$ structures on a manifold and see that they produce K-theory orientations on 
the tangent bundle. Given a Riemannian manifold $(X,g)$, we can form the (complex) {\bf Clifford bundle}\footnote{Of course we also have the real Clifford
bundle $\RCl(X)\to X$, but we shall not need it for our applications.  More generally, we can define Clifford bundles $\RCl(V) \to X$ and $\Cl(V)\to X$ whenever
$V \to X$ is a vector bundle with inner product.}
$$
	\Cl(X) \to X, \quad \Cl(X)_x := \Cl(T_x X, g_x), \text{ for all $x \in X$,}
$$
which is a bundle of complex Clifford algebras of dimension $n = \dim(X)$.  A complex vector bundle $E \to X$ is called a {\bf Clifford module} if it carries a fiberwise
action
$$
	\cl : \Cl(X) \to \End(E).
$$
Such an action, if it exists, will be called {\bf Clifford multiplication}.

Of course, over each $x \in X$, any Clifford module decomposes as a direct sum of irreducible modules over $\Cl(X)_x$, but this is not necessarily true globally.
This leads us to the notion of spin and spin$^c$ structures on $X$.  

Let $P_{\SO}(X)\to X$ be the frame bundle of $X$, i.e.\ the principal $\SO_n$ bundle to which $TX$ is associated:
$$
	TX = P_{\SO}(X)\times_{\SO_n} \R^n.
$$
$X$ is called a {\bf spin manifold} if there exists a principal $\Spin_n$ bundle $P_{\Spin}(X)$ and a bundle map
$$
\xymatrix{ \Spin_n \ar[r]\ar[d] &\SO_n\ar[d] \\
	P_{\Spin}(X) \ar[r]\ar[dr] &P_{\SO}(X)\ar[d] \\
	& X
}
$$
which is a 2-sheeted cover of $P_\SO(X)$.  $P_\Spin(X)$ is called a {\bf spin structure on $X$}.\footnote{Similarly, a general vector bundle with inner 
product $V \to X$ admits a spin structure whenever $P_\SO(V)$ admits a 2-sheeted cover $P_\Spin(V)$.}  The obstruction to obtaining such a cover of $P_\SO(X)$ is the 
second Stiefel-Whitney class\footnote{This is straightforward to see by trying to patch $P_\Spin(X)$ together over a trivializing cover.  In order to do so,
we must have a Cech ``cohomology class'' (I'm using quotes since the coefficients are in a nonabelian group; nevertheless $H^1(X; \SO_n)$ is a based set)
in $H^1(X; \SO_n)$ which is the image of a class in $H^1(X; \Spin_n)$.  Using the long exact sequence associated to 
$$
	1 \to \Z_2 \to \Spin_n \to \SO_n \to 1,
$$
the image of this class in $H^2(X; \Z_2)$ is exactly $w_2(X) \in H^2(X; \Z_2)$.} of $X$:
$$
	X \text{ is spin iff } w_2(X) \equiv 0,
$$
and if $X$ is spin, the possible spin structures of $X$ are parametrized by $H^1(X, \Z_2)$.

A {\bf spin$^c$ structure on $X$} consists of a complex line bundle $L \to X$ and a lift
$$
\xymatrix{ \Spin^c_n \ar[r]\ar[d] &\SO_n\times \mathrm{U}_1\ar[d] \\
	P_{\Spin^c}(X,L) \ar[r]\ar[dr] &P_{\SO}(X)\times P_{U_1}(L)\ar[d] \\
	& X
}
$$
which is a 2-sheeted covering of $P_\SO(X)\times P_{U_1}(L)$, where $P_{U_1}(L)$ is the structure bundle of $L$.  $X$ is called a {\bf spin$^c$ manifold}
if such a lift exists.  The obstruction to obtaining $P_{\Spin^c}(X,L)$ is the class $w_2(X) + c_1(L)(\text{mod }2)$; thus
$$
	X \text{ is spin$^c$ iff } w_2(X) = \alpha (\text{mod }2) \text{ for some }\alpha \in H^2(X, \Z)
$$
Being spin$^c$ is a weaker condition than being spin: 
\begin{prop}
If $X$ is spin, then it has a canonical spin$^c$ structure associated to the trivial line bundle, so 
$$
	X \text{ spin } \implies X \text{ spin$^c$.}
$$
\end{prop}
Additionally, any (almost) complex manifold has a canonical spin$^c$ structure. 
\begin{prop}
If $X$ is an almost complex manifold, then $w_2(X) = c_1(X) (\text{mod }2)$, and $X$ has a canonical spin$^c$ structure
associated to the determinant line bundle $\swedge^n_{\C} TX$ (which satisfies $c_1\parens{\swedge^nTX} = c_1(X)$).
\end{prop}
Note that if $X$ is both spin and almost complex, the spin$^c$ structure coming from the spin structure is generally {\em not} the same as
the one coming from the almost complex structure.

The importance of spin$^c$ structures is the following proposition, which says that, given a spin$^c$ structure, Clifford modules are globally reducible, 
and in bijection with the set of $\Cl_n$ modules.
\begin{prop}
If $X$ is spin$^c$, then every Clifford module $E \to X$ has the form
$$
	E = P_{\Spin^c}(X,L) \times_{\sigma} F,
$$
where $\sigma : \Spin^c_n \to \GL(F)$ is a representation of $\Spin^c_n$ which extends to a representation of $\Cl_n$ (here $n = \dim(X)$).  
\end{prop}
Turning the construction around, we obtain Clifford modules over a spin$^c$ manifold $X$ for every representation of $\Cl_n$; in particular for $\dim(X) = 2n$, 
we have the complex {\bf spinor bundle}
$$
	\bbS(X) = \bbS^+(X)\oplus\bbS^-(X)  = P_{\Spin^c}(X,L) \times_{\rho^+_{1/2}\oplus \rho^-_{1/2}} \bbS^+\oplus \bbS^-
$$
with the (graded) action
$$
	\cl : \Cl(X) \to \End_{\mathrm{gr}}(\bbS^+(X)\oplus\bbS^-(X)).
$$

Finally, we can get to the main point about spin$^c$ structures, which is that they allow us to construct K-theory orientation classes\footnote{In fact the proposition
is valid for general vector bundles $V \to X$ with spin$^c$ structure; the analogous element $\mu = [\pi^\ast\bbS^+(V), \pi^\ast \bbS^-(V), \cl] \in K^0_c(V)$ 
is a Thom class.}
 for $T^\ast X \to X$.
\begin{prop}
If an even dimensional manifold $X$ has a spin$^c$ structure and $\bbS(X) = \bbS^+(X)\oplus\bbS^-(X)$ is the bundle of spinors, then
$$
	\mu = [\pi^\ast \bbS^+(X),\pi^\ast \bbS^-(X), \cl] \in K^0_c(T^\ast X)
$$
is an orientation/Thom class for complex K-theory, so $K^\ast_c(T^\ast X)$ is freely generated by $\mu$ as a module over $K^\ast(X)$, and
$$
	K^\ast_c(T^\ast X) \cong K^\ast(X).
$$
\end{prop}
Note that $T^\ast X \subset \Cl^1(X)$, so $\cl(\xi) : \bbS^\pm(X)_x \to \bbS^\mp(X)_x$ for $(x,\xi) \in T^\ast X$; moreover, this multiplication is invertible 
with inverse $\abs{\xi}^{-2}\cl(\xi)$ provided $\xi \neq 0$.

As a side remark, let me point out that an analogous theorem is true for spin structures and real K-theory: If $X$ is spin and $8n$ dimensional, then 
$[\pi^\ast \bbS^+(X), \pi^\ast\bbS^-(X), \cl] \in KO^0_c(T^\ast X)$ is an orientation class, and $KO^\ast(X) \cong KO_c^\ast(T^\ast X)$.

Note that if $X$ is almost complex, with the corresponding spin$^c$ structure, then we can identify $\bbS(X)$ with $\swedge^\ast_\C T^\ast X$; in this case,
$$
	\bbS^\pm(X) \cong \swedge^\mathrm{even/odd}_\C T^\ast X
$$
and $\cl(\xi) = \xi \smwedge \cdot - \xi^\ast \iprod \cdot$ under this identification.  Thus we recover the Thom element
$$
	\mu = [\pi^\ast \swedge^\mathrm{even} V, \pi^\ast \swedge^\mathrm{odd} V, \cl] \in K^0_c(V)
$$
for complex bundles, and we see that the Thom isomorphism for such bundles can be thought of as a special case of the isomorphism for spin$^c$ bundles.

Finally, let us briefly discuss what this looks like in the setting of a fibration $X \to Z$.  In this case the relevant Clifford bundle is
$$
	\Cl(X/Z) = \Cl(T(X/Z), g) \to X,
$$
and the fibration is oriented as long as $T(X/Z) \to X$ admits a spin$^c$ structure.  Indeed, if it does, we have the orientation class
$$
	\mu = [\pi^\ast \bbS^+(X/Z), \pi^\ast \bbS^-(X/Z), \cl] \in K^0_c(T^\ast (X/Z))
$$
constructed from the bundles of spinors $\bbS^\pm(X/Z) \to X$.

\section{Dirac operators} \label{S:dirac}
In this section we'll see that the orientation classes discussed above are in fact the symbols of particularly nice (families of) elliptic differential operators.
Let $E \to X$ be any Clifford module over $X$.  Suppose $E$ is endowed with a connection $\nabla : C^\infty(X; E) \to C^\infty(X; T^\ast X \otimes E)$ such that
$$
	\nabla(\cl(\xi) s) = \cl(\nabla^{\mathrm{LC}} \xi) s + \cl(\xi) \nabla s
$$
where $\nabla^{\mathrm{LC}}$ is the Levi-Civita connection, which extends canonically to $\Cl(X) \to X$.  Such a connection (which always exists) is called 
a {\bf Clifford connection} on $E$.  

Given such data, we can construct a canonical first order, elliptic, differential operator $D \in \Diff^1(X; E)$ called 
a {\bf Dirac operator}; at a point $x \in X$, $D$ is defined by
$$
	D_p = \sum_i \cl(e_i) \nabla_{e_i}, \quad \text{ for an orthonormal basis $\set{e_i}$ of $T_xX$.}
$$
\begin{prop}
Such a Dirac operator is an elliptic, essentially self adjoint (on $L^2(X; E)$) operator with principal symbol
$$
	\sigma(D)(\xi) = i\cl(\xi).
$$
If $E = E^+\oplus E^-$ is a graded $\Cl(X)$ module, then $D$ has the form
$$
	D = \begin{pmatrix} 0 & D^- \\ D^+ & 0\end{pmatrix}
$$
with $D^+$ and $D^-$ mutual adjoints.
\end{prop}
Note that $\sigma(D^2) = \sigma(D)^2 = \abs{\xi}^2\id$, so $D^2$ is a Laplacian operator on $E$.

If $X$ is a spin$^c$ manifold, we can form a canonical {\bf spin$^c$ Dirac operator} 
$$
	\dirD = \begin{pmatrix} 0 & \dirD^- \\ \dirD^+ & 0\end{pmatrix}  \in \Diff^1(X; \bbS^+(X)\oplus \bbS^-(X))
$$ 
acting on the spinors $\bbS(X)$, and it follows that
$$
	[\pi^\ast \bbS^+(X), \pi^\ast \bbS^-(X), \sigma(\dirD^+)] \in K^0_c(T^\ast X)
$$
is an orientation class for K-theory.  This realizes the Thom isomorphism as follows.  We can always twist $\dirD$ by a vector bundle $E \to X$ by trivially extending
the Clifford action to the bundle $\bbS(X)\otimes E$ and taking a product connection to get $\dirD_E \in \Diff^1(X; \bbS(X)\otimes E)$.  Then for an element
$$
	[E] - [F] \in K^0(X),
$$
the image under the Thom isomorphism $K^0(X) \to K^0_c(T^\ast X)$ is the element\footnote{The better way to write this is to use $\Z_2$ gradings everywhere.  
Let $\bbE = E\oplus F$ considered as a graded vector bundle and form $\bbS(X)\hat\otimes \bbE$, where $\hat\otimes$ denotes the graded tensor product.  
Then $\dirD$ extends to a graded, twisted operator $\dirD_\bbE$, and the orientation class is given by 
$[\pi^\ast (\bbS(X)\hat\otimes \bbE)^+, \pi^\ast (\bbS(X)\hat\otimes \bbE)^-, \sigma(\dirD_\bbE^+)] \in K^0_c(T^\ast X)$.  We'll talk more about gradings
in section \ref{S:higher}.}
$$
	[\pi^\ast \bbS^+(X)\otimes E, \pi^\ast\bbS^-(X)\otimes E, \sigma(\dirD_E)] - [\pi^\ast \bbS^+(X)\otimes F, \pi^\ast\bbS^-(X)\otimes F, \sigma(\dirD_F)] 
		\in K^0_c(T^\ast X).
$$

For a complex manifold $X$, you have probably already met the canonical spin$^c$ Dirac operator.  Indeed, using the identifications 
$\bbS(X) \cong \swedge^\ast_\C T^\ast X \cong \swedge^{0,\ast} T^\ast X$, one can easily see that 
$$
	\dirD^+ = \overline{\pa} + \overline{\pa}^\ast \in \Diff^1(X; \swedge^{0,\mathrm{even}}T^\ast X, \swedge^{0,\mathrm{odd}}T^\ast X)
$$
is just the {\bf Dolbeault operator} acting from even to odd harmonic forms.

Finally, in the case of an oriented fibration $X \to Z$, we construct in precisely the same way the canonical {\em family of spin$^c$ Dirac operators}
$$
	\dirD \in \Diff^1(X/Z; \bbS^+(X/Z)\oplus \bbS^-(X/Z))
$$
and of course
$$
	[\pi^\ast \bbS^+(X/Z), \pi^\ast \bbS^-(X/Z), \sigma(\dirD)] \in K^0_c(T^\ast(X/Z))
$$
is the Thom class.  

This gives the analytical realization of the Gysin map $K^0(X) \to K^0(Z)$; namely, it coincides with the analytical index of the family of spin$^c$ Dirac operators,
twisted by the given element in $K^0(X)$:
$$
	K^0(X) \ni [E] - [F] \mapsto \ind( \dirD_E - \dirD_F) \in K^0(Z).
$$

\section{Higher Index} \label{S:higher}
We will develop two pictures of the higher $K$-groups of a manifold $X$, in analogy to the two we've developed for $K^0(X)$, namely, the Grothendieck group
of vector bundles, and the classifying space consisting of Fredholm operators.  Really this whole story is a bit more interesting in the case of real K-theory,
and much of what we describe below will be valid if one replaces $K^\ast(X)$ by $KO^\ast(X)$ and $\Cl_\ast$ by $\RCl_\ast$ (with the obvious exception of 2 periodicity
of $\Cl$ modules, which would be replaced by an analogous 8-fold periodicity of $\RCl$ modules).

Fix $k$ for a moment, and consider the semigroup of $\Cl_k$ modules.  Of course this can be completed to a group by the usual Grothendieck construction, and
we denote the {\bf Grothendieck group of $\Cl_k$ modules} by $\cM_k$.  Now, the inclusion $i : \R^k \hookrightarrow \R^{k+1}$ induces an injective algebra 
homomorphism $i : \Cl_k \hookrightarrow \Cl_{k+1}$, which gives a restriction operation
$$
	i^\ast : \cM_{k+1} \to \cM_k
$$
on Clifford modules.  It turns out that the interesting object to consider is $\cM_k/i^\ast\cM_{k+1}$.

Actually, it is more convenient at this point to work in terms of graded modules.  Thus, let $\hcM_k$ denote the {\bf Grothendieck group of graded $\Cl_k$ modules}.  
Again we have a restriction
$$
	i^\ast : \hcM_{k+1} \to \hcM_k,
$$
and from the equivalence between graded $\Cl_k$ modules and ungraded $\Cl_{k-1}$ modules, we have
$$
	\hcM_k / i^\ast \hcM_{k+1} \cong \cM_{k-1} / i^\ast \cM_k. 
$$
Furthermore, it is easy to check using the representation theory of complex Clifford algebras, that we have the following periodicity (related of course to Bott
periodicity)
$$	
	\hcM_k / i^\ast \hcM_{k+1} \cong \cM_{k-1} / i^\ast \cM_k = \begin{cases} \Z & \text{if $k$ is even} \\ 0 & \text{if $k$ is odd.} \end{cases}
$$

Let $W = W^0 \oplus W^1 \in \cM_k$, and form the trivial bundles $E^i = D^k \times W^i$ over the unit disk $D^k \subset \R^k$.  We can form the element
$$
	\set{E^0,E^1,\cl(\cdot)} \in K^0(D^k, S^k)
$$
where $\cl(\cdot) : S^k \subset \R^k\setminus \set{0} \subset \Cl_k \to \Iso(W^0,W^1)$.  This bundle isomorphism over $S^k$ can be shown to extend over $D^k$ if and 
only if $W$ actually comes from a $\Cl_{k+1}$ module.  Thus one obtains the celebrated result of Atiyah, Bott and Shapiro \cite{atiyah1964clifford}.
\begin{thm}[Atiyah-Bott-Shapiro] The above construction gives a graded ring isomorphism
$$
	\hcM_\ast/i^\ast\hcM_{\ast+1} \cong K^0(D^\ast,S^\ast) = K^{-\ast}(\pt).
$$
\end{thm}
\cite{atiyah1964clifford} contains an analogous result for real K-theory and $\RCl$ modules.  
While the above looks like an appealing way to prove Bott periodicity from the more obvious periodicity of Clifford modules, it is not actually so.  Indeed, 
the ABS result uses periodicity of K-theory in the proof.

This leads to the analogue of the vector bundle representation of $K^0(X)$.  Namely, elements of $K^k(X)$ can be represented\footnote{I'm not sure of a good reference
for this explicit representation of higher K-theory.  It is implicit in Karoubi's formulation, but he takes the algebraic approach, with projective $C^0(X)$ modules
instead of vector bundles.} as (isomorphism classes of) bundles of graded $\Cl_k$ modules, modulo those which admit a graded $\Cl_{k+1}$ action:
$$
	K^k(X) = \set{V^0\oplus V^1 \to X \;;\; \Cl_k \to \End_{\mathrm{gr}}(V^0\oplus V^1)}/\set{\Cl_{k+1} \to \End_{\mathrm{gr}}(V^0\oplus V^1)}
$$
It is an instructive exercise to recover the vector bundle representation of $K^0(X)$ from this definition.  Indeed, since $\Cl_0 = \C$ with the trivial grading and 
$\Cl_1 = \C\oplus\C$, we see that graded $\Cl_0$ modules are just vector bundles of the form $E\oplus F$, which extend to $\Cl_1$ modules only if
$E \cong F$ (since the action of the generator of the 1-graded part of $\Cl_1$ must be an isomorphism: $0\oplus 1 : E\oplus F \stackrel{\cong}{\to} F\oplus E$).
Thus we have the equation $E\oplus E = 0 \in K^0(X)$ or equivalently $- E\oplus 0 = 0\oplus E$, and we see that we can identify $[E] - [F]$ in the old representation
with $E\oplus F$ in this new representation.

Next we can generalize the Fredholm operator representation of $K^0(X)$, following Atiyah and Singer's paper \cite{atiyah1969index}.  Let $H = H^0\oplus H^1$ be a
graded, separable, infinite dimensional Hilbert space, and assume $H$ is a module over $\Cl_k$ (for a given $k$, but then the action can be extended for all $k$).
Let
$$
	\Fred_k(H) = \set{P \in \Hom(H^i; H^{i+1})\;;\; P \text{ Fredholm and commutes with $\Cl_k$}}
$$
be the space of graded (i.e.\ acting as 1-graded elements) Fredholm operators\footnote{Actually this is a bit of an oversimplification.  When $k$ is odd, $\Fred_k(H)$
can be separated into open components $\Fred_k^+$, $\Fred_k^-$ and $\Fred_k^\ast$ consisting of operators which are essentially positive (meaning positive off of a 
finite dimensional subspace), essentially negative, or neither.  The first two are contractible, and we take $\Fred_k^\ast(H)$ in this case.} 
commuting (in the graded sense) with $\Cl_k$.  We will call these the {\bf $\Cl_k$-linear Fredholm operators}.  These form a classifying space for $K^{-k}(\ast)$:
\begin{thm}[Atiyah-Singer]
There is an explicit homotopy equivalence
$$
	\Fred_k(H) \simeq \Omega \Fred_{k-1}(H)
$$
for all $k$, and therefore
$$
	[X, \Fred_k(H)] = K^{-k}(X)
$$
\end{thm}
Again, for $P \in [X, \Fred_k(H)]$ one can morally take $[\ker P] \in K^{-k}(X)$, since at each point $\ker P$ is a graded $\Cl_k$ module.  A stabilization procedure
would be required to make this precise, and I have to admit I've never seen it written down, though I'm sure it's possible.

Finally, one can make the Atiyah-Singer index construction go through in this case (again, I've not seen this written explicitly, but reliable sources assure me it's 
true!).  Namely, given a fibration $X \to Z$, if one has a {\bf $\Cl_k$ linear family of elliptic differential operators} 
$$
	P \in \Diff^l(X/Z; E^0\oplus E^1),
$$
meaning that $P = \begin{pmatrix} 0 & P_1 \\ P_0 & 0\end{pmatrix}$ is graded and commutes in the graded sense with an action 
$\Cl_k \to \End_{\mathrm{gr}}(E^0\oplus E^1)$, then
$$
	\ind(P) = [\ker P_0 \oplus \ker P_1] \in K^k(Z)
$$
and that this index coincides with the Gysin map
$$
	\ind = p_! : K^k_c(T^\ast X/Z) \to K^k(Z).
$$

Of course, since $K^1(\pt) = 0$, operators on a manifold $X$ never have any interesting odd index (however a family of operators might, provided $K^1(Z) \neq 0$).
For real K-theory, however, this can be an interesting and useful concept.  For instance, the Kervaire semicharacteristic on a ($4k +1$) manifold $X$ can be 
computed (see \cite{lawson1989spin}) as the odd index in $KO^1(\pt) = \Z_2$ of a $\RCl_1$ linear elliptic differential operator on $X$!\footnote{Sorry about all 
the footnotes.}

\bibliographystyle{amsalpha}
\bibliography{references}

\providecommand{\bysame}{\leavevmode\hbox to3em{\hrulefill}\thinspace}
\providecommand{\MR}{\relax\ifhmode\unskip\space\fi MR }
\providecommand{\MRhref}[2]{%
  \href{http://www.ams.org/mathscinet-getitem?mr=#1}{#2}
}
\providecommand{\href}[2]{#2}
\begin{thebibliography}{ABS64}

\bibitem[ABS64]{atiyah1964clifford}
M.F. Atiyah, R.~Bott, and A.~Shapiro, \emph{{Clifford modules}}, Topology
  \textbf{3} (1964), no.~3, 38.

\bibitem[AS68]{atiyah1968index_I_III}
M.F. Atiyah and I.~Singer, \emph{{Index theorem of elliptic operators, I,
  III}}, Annals of Mathematics \textbf{87} (1968), 484--530.

\bibitem[AS69]{atiyah1969index}
\bysame, \emph{{Index theory for skew-adjoint Fredholm operators}},
  Publications Math{\'e}matiques de l'IH{\'E}S \textbf{37} (1969), no.~1,
  5--26.

\bibitem[AS71]{atiyah1971index_IV}
\bysame, \emph{{The index of elliptic operators: IV}}, Annals of Mathematics
  (1971), 119--138.

\bibitem[Ati67]{atiyah1967ktheory}
M.F. Atiyah, \emph{{K-theory: lectures by MF Atiyah; Notes by DW Anderson}}.

\bibitem[LM89]{lawson1989spin}
H.B. Lawson and M.L. Michelsohn, \emph{{Spin geometry}}, Princeton University
  Press, 1989.

\end{thebibliography}

\end{document}